\documentclass[a4paper,10pt]{article}
\usepackage[utf8x]{inputenc}
\usepackage{tracefnt,amsmath,tabu,array}
\usepackage{amssymb,graphicx,setspace,amsfonts,amsbsy}
\usepackage{pifont,latexsym,ifthen,amsthm,rotating,calc,textcase,booktabs,cancel,slashed}

\newtheorem{theorem}{Theorem}[section]
\newtheorem{lemma}[theorem]{Lemma}

\newtheorem{proposition}[theorem]{Proposition}

\newcommand{\filledbox}{\leavevmode
  \hbox to.77778em{%
  \hfil\vbox to.675em{\hrule width.6em height.6em}\hfil}}

\newcommand{\Rm}{{\mathbb R}}

\begin{document}
\tabulinesep=1.0mm
\title{Asymptotic behaviour of non-radiative solution to the wave equations}

\author{Liang Li, Ruipeng Shen, Chenhui Wang and Lijuan Wei\\
Centre for Applied Mathematics\\
Tianjin University\\
Tianjin, China
}

\maketitle

\begin{abstract}
 In this work we consider weakly non-radiative solutions to both linear and non-linear wave equations. We first characterize all weakly non-radiative free waves, without the radial assumption. Then in dimension 3 we show that the initial data of non-radiative solutions to a wide range of nonlinear wave equations are similar to those of non-radiative free waves in term of asymptotic behaviour. 
\end{abstract}

\section{Introduction}

\subsection{Background and main topics}

\paragraph{Channel of energy} The channel of energy method plays an important role in the study of asymptotic behaviour of solutions to non-linear wave equations in the past decade. This method mainly discusses the distribution of energy as time tends to infinity. More precisely, if $u$ is a solution to either linear or non-linear wave equation defined for all time, then the following limits are considered for a given constant $R$.
\[ 
 \lim_{t\rightarrow \pm \infty} \int_{|x|>R+|t|} |\nabla_{t,x} u(x,t)|^2 dx. 
\] 
Here for convenience we use the notation $\nabla_{t,x} u = (u_t, \nabla u)$. This theory was first established for solutions to homogeneous linear wave equation, i.e. free waves, then applied to the study of non-linear wave equations. Please see, for instance,  C\^{o}te-Kenig-Schlag \cite{channeleven}, Duyckaerts-Kenig-Merle \cite{tkm1, oddtool} and Kenig-Lawrie-Schlag \cite{channel5d} for linear theory; and Duyckaerts-Kenig-Merle \cite{se, oddhigh} for the applications of the channel of energy on soliton resolution of focusing wave equation. 

\paragraph{Non-radiative solutions} A crucial part of the channel of energy theory is to discuss the property of non-radiative solutions. Let $u$ be a solution to the wave equation with a finite energy. We call it a non-radiative solution if and only if 
\[ 
 \lim_{t\rightarrow \pm \infty} \int_{|x|>|t|} |\nabla_{t,x} u(x,t)|^2 dx = 0. 
\] 
We may also consider a more general case. We call a solution $u$ to be $R$-weakly non-radiative if and only if 
\[ 
 \lim_{t\rightarrow \pm \infty} \int_{|x|>R+|t|} |\nabla_{t,x} u(x,t)|^2 dx = 0. 
\] 
Let us first consider (weakly) non-radiative solutions to the homogeneous linear wave equation in $\Rm^d$. It has been proved that any non-radiative free wave must be zero, see Duyckaerts-Kenig-Merle \cite{dkmnonradial, oddtool}.  All radial weakly non-radiative free waves have also been well understood. The following result was first proved for odd dimensions $d\geq 3$ by Kenig et al \cite{channel} then generalized to the even dimensions $d\geq 2$ in Li-Shen-Wei 
\cite{shenradiation}. 
 
\begin{proposition} [Radial weakly non-radiative solutions] \label{radial structure}
 Let $d\geq 2$ be an integer and $R>0$ be a constant. If initial data $(u_0,u_1) \in \dot{H}^1 \times L^2$ are radial, then the corresponding solution to the homogeneous linear wave equation $u$ is $R$-weakly non-radiative, i.e. 
 \[
  \lim_{t\rightarrow \pm \infty} \int_{|x|>|t|+R} |\nabla_{t,x} u(x,t)|^2 dx = 0,
 \]
if and only if the restriction of $(u_0,u_1)$ in the region $\{x\in \Rm^d: |x|>R\}$ is contained in
\[
  \hbox{Span} \left\{(r^{2k_1-d},0), (0,r^{2k_2-d}): 1\leq k_1 \leq \left\lfloor \frac{d+1}{4} \right\rfloor, 1\leq k_2 \leq \left\lfloor \frac{d-1}{4} \right\rfloor\right\}
\]
Here the notation $\lfloor q\rfloor$ is the integer part of $q$. In particular, all radial $R$-weakly non-radiative solution in dimension $2$ are supported in $\{(x,t): |x|\leq |t|+R\}$. 
\end{proposition}

\paragraph{Goals of this work} The aim of this paper is two-fold. The first goal of this paper is to characterize all (possibly non-radial) initial data so that the corresponding solutions to free wave equation are $R$-weakly non-radiative. For convenience we define 
\[
 P(R) \doteq \left\{(u_0,u_1)\in \dot{H}^1 \times L^2(\Rm^d): \lim_{t\rightarrow \pm \infty} \int_{|x|>R+|t|} |\nabla_{t,x} \mathbf{S}_L (u_0,u_1)|^2 dx = 0\right\}.
\]
Here $\mathbf{S}_L(u_0,u_1)$ is the corresponding solution of the free wave equation with given initial data $(u_0,u_1)$. We will give a decomposition of every element $(u_0,u_1)\in P(R)$ in term of spherical harmonic functions, whose details are given in Section \ref{sec: non radial space PR}. The second goal is to show that in the 3-dimensional case any weakly non-radiative solution to a wide range of non-linear wave equations share the same asymptotic behaviour as weakly non-radiative free waves, as given in Section \ref{sec: nonlinear case}. Our argument depends on a suitable decay estimate of weakly non-radiative free waves in the exterior region $\{(x,t): |x|>|t|+R\}$. The decay estimates of this kind are clearly true for radial non-radiative solutions, as given in Proposition \ref{radial structure}. Although we expect that a similar estimate holds for non-radial non-radiative solutions in all dimensions $d \geq 2$ as well, this has been proved only in dimension 3, as far as the author knows. This is why we have to restrict our discussion to dimension 3 in this work. 

\section{The characteristics of $P(R)$} \label{sec: non radial space PR}

In this section we give an explicit expression of the element in the space $P(R)$. We use spherical harmonics and follow a similar argument as given in Duychaerts-Kenig-Merle \cite{exteriorW}. Let us first give a brief review on some basic properties of spherical harmonics. We recall that the eigenfunctions of the Laplace-Beltrami operator on $\mathbb{S}^{d-1}$ are exactly the homogeneous harmonic polynomials of the variables $x_1, x_2, \cdots, x_d$. Such a polynomial $\Phi$ of degree $\nu$ satisfies 
\[
 -\Delta_{\mathbb{S}^{d-1}} \Phi = \nu(\nu+d-2) \Phi. 
\]
We choose a Hilbert basis $\{\Phi_k(\theta)\}_{k\geq 0}$ of the operator $-\Delta_{\mathbf{S}^{d-1}}$ on the sphere $\mathbb{S}^{d-1}$. Here we assume that the harmonic polynomial $\Phi_k$ is of degree $\nu_k$. In particular we assume $\nu_0 = 0$ and $\nu_k > 0$ if $k \geq 1$. Next we give the statement of our first main result. We start by the odd dimensional case and then deal with the even dimensional case. Please note that a similar result for odd dimensions has been proved in C\^{o}te-Laurent \cite{newradiation} by the Radon transform. The novelty of our result includes
\begin{itemize}
 \item We give an $L^2$ decay estimate of $\partial_r u_0$ near infinity in addition;
 \item The argument works for even dimensions as well, with minor modifications. 
\end{itemize}

\subsection{Odd dimensions}

\begin{proposition} \label{structure of PR}
 Assume that $d\geq 3$ is an odd integer and $\mu = (d-1)/2$. Then $(u_0,u_1) \in P(R)$ is equivalent to saying that there exist two sequences of polynomials $\{P_k(z)\}_{k\geq 0}$ and $\{Q_k(z)\}_{k\geq 0}$ of the following form ($A_{k,k_1}, B_{k,k_2}$ are constants)
 \begin{align*}
 &P_k(z) = \sum_{1\leq k_1 \leq \lfloor\frac{\mu+\nu_k+1}{2}\rfloor} A_{k,k_1} z^{\mu+1+\nu_k-2k_1};& &Q_k(z) = \sum_{1\leq k_2 \leq \lfloor\frac{\mu+\nu_k}{2}\rfloor} B_{k,k_2} z^{\mu+\nu_k-2k_2};&
\end{align*}
with
\begin{align*}
 &\sum_{k=0}^\infty \int_0^{1/R} \left( \nu_k (d-2+\nu_k) \left|P_k(z)\right|^2 + |zP'_k(z)|^2 \right) dz < +\infty;& &\sum_{k=0}^\infty \int_0^{1/R} \left|Q_k(z)\right|^2 dz < +\infty;&
\end{align*}
so that 
 \begin{align}
 &u_0(r,\theta) = \sum_{k=0}^\infty r^{-\mu} P_k(1/r) \Phi_k(\theta), \quad r>R;& &u_1(r,\theta) = \sum_{k=0}^\infty r^{-\mu-1} Q_k(1/r) \Phi_k(\theta), \quad r>R.& \label{bh map}
\end{align}
Here the first identity holds for every fixed $r>R$ in the sense of $L^2(\mathbb{S}^{d-1})$ convergence. The second one holds in the sense of $L^2(\{x: |x|>R\})$ convergence. In addition, we have 
\begin{itemize}
 \item[(i)] The derivative of $u_0$ can be given by 
 \begin{align*}
  \nabla_x u_0 (r,\theta) = \sum_{k=0}^\infty r^{-\mu-1} \left\{P_k(1/r) \nabla_\theta \Phi_k(\theta) - [\mu P_k(1/r) + (1/r) P'_k(1/r)] \Phi_k(\theta) \theta \right\}.
 \end{align*}
This identities holds in the sense of $L^2(\{x: |x|>R\})$ convergence. Here we naturally embed $\nabla_\theta \Phi_k (\theta)$ into $\Rm^d$ by the identity $\nabla_\theta \Phi_k(\theta) = \nabla_x \Phi_k(\theta)$. In the right hand side of this identity we understand $\Phi_k$ as a function defined in $\Rm^d\setminus \{0\}$ by polar coordinates. 
 \item[(ii)] The norms of $(u_0,u_1)$ can be determined by $P_k(z)$ and $Q_k(z)$'s:
 \begin{align*}
  \|\slashed{\nabla} u_0\|_{L^2(\{x: |x|>R\})}^2 & = \sum_{k=1}^\infty \nu_k (d-2+\nu_k) \int_0^{1/R} \left|P_k(z)\right|^2 dz;\\
  \|u_1\|_{L^2(\{x: |x|>R\})}^2 & =  \sum_{k=0}^\infty \int_0^{1/R} \left|Q_k(z)\right|^2 dz; \\
  \|\partial_r u_0\|_{L^2(\{x: |x|>R\})} & =  \sum_{k=0}^\infty \int_0^{1/R} |zP'_k(z)+\mu P_k(z)|^2 dz < +\infty.
 \end{align*} 
 \item[(iii)] The derivative $\partial_r u_0$ satisfies the following decay estimates ($R_1 \geq 2R$)
 \begin{align*}
  \int_{|x|>R_1} |\partial_r u_0(x)|^2 dx & \lesssim (R/R_1) \int_{|x|>R} |\nabla u_0(x)|^2  dx;\\
  \int_{|x|>R_1} |\partial_r u_0^\ast (x)|^2 dx  & \lesssim (R/R_1) \int_{|x|>R} |\slashed{\nabla} u_0(x)|^2 dx. 
 \end{align*}
 Here $u_0^\ast$ is the non-radial part of $u_0$ defined by $u_0^\ast = u_0 - r^{-\mu} P_0(1/r) \Phi_0$. 
\end{itemize}
\end{proposition}

\paragraph{Proof} The rest of this subsection is devoted to the proof of this proposition. The proof consists of three parts: Step one, we first show that any element in $P(R)$ can be written as in \eqref{bh map}. Step two, we show any initial data given by \eqref{bh map} is indeed contained in $P(R)$. Finally in Step three we prove the identities and inequalities in the proposition. 

\paragraph{Step one} Let us consider 
\[
 u_k(r,t) = r^{-\nu_k} \int_{\mathbb{S}^{d-1}} u(r\theta, t) \Phi_k(\theta) d\theta.
\]
Let $ \Box = \partial_t^2 - \partial_r^2 - \frac{d+2\nu_k-1}{r} \partial_r$. A basic calculation shows 
\begin{align*}
 \Box u_k & = (\Box r^{-\nu_k}) \int_{\mathbb{S}^{d-1}} u(r\theta, t) \Phi_k(\theta) d\theta + r^{-\nu_k} \int_{\mathbb{S}^{d-1}} \Box u(r\theta, t) \Phi_k(\theta) d\theta\\
 & \qquad - 2 \partial_r(r^{-\nu_k}) \int_{\mathbb{S}^{d-1}} \partial_r u(r\theta,t) \Phi_k(\theta) d\theta\\
 & = r^{-\nu_k} \int_{\mathbb{S}^{d-1}} \left(\partial_t^2 - \partial_r^2 - \frac{d-1}{r}\partial_r \right) u (r\theta, t) \Phi_k(\theta) d\theta + \nu_k(d-2+\nu_k) r^{-2} u_k\\
 & = r^{-\nu_k} \int_{\mathbb{S}^{d-1}} r^{-2} \Delta_{\mathbf{S}^{d-1}} u (r\theta, t) \Phi_k(\theta) d\theta + \nu_k(d-2+\nu_k) r^{-2} u_k\\
 & = r^{-\nu_k-2} \int_{\mathbb{S}^{d-1}}  u (r\theta, t) \Delta_{\mathbf{S}^{d-1}} \Phi_k(\theta) d\theta + \nu_k(d-2+\nu_k) r^{-2} u_k\\
 & = 0.
\end{align*}
Thus if $u_k$ is viewed as a radial function defined on $\Rm^{d+2\nu_k}$, it satisfies the free wave equation 
\[ 
\partial_t^2 u_k - \Delta_{\Rm^{d+2\nu_k}} u_k = 0, \qquad |x|>0.
\]
In addition, we have 
\begin{align*}
 \int_{R+|t|}^{\infty} & \left(|\partial_r u_k(r,t)|^2 + |\partial_t u_k(r,t)|^2 \right) r^{d+2\nu_k-1} dr\\
& = \int_{R+|t|}^\infty \left(\left|\int_{\mathbb{S}^{d-1}} \partial_r u(r\theta,t) \Phi_k(\theta) d\theta \right|^2 +\left|\int_{\mathbb{S}^{d-1}} \partial_t u(r\theta,t) \Phi_k(\theta)d\theta \right|^2 \right) r^{d-1} dr\\
& \leq  \int_{R+|t|}^\infty \left(\int_{\mathbb{S}^{d-1}} |\partial_r u(r\theta,t)|^2 d\theta + \int_{\mathbb{S}^{d-1}} |\partial_t u(r\theta,t)|^2 d\theta \right) r^{d-1} dr\\
& \lesssim \int_{|x|>R+|t|} |\nabla_{x,t} u(x,t)|^2 dx. 
\end{align*}
Thus $u_k$ is also a weakly non-radiative solution. According to the explicit expression of radial non-radiative solutions, there exist constants $A_{k, k_1}$ and $B_{k, k_2}$, so that 
\begin{align*}
 u_k(r,0) & = \sum_{1\leq k_1 \leq \frac{\mu+\nu_k+1}{2}} A_{k,k_1} r^{-d-2\nu_k+2k_1} = r^{-\mu-\nu_k} P_k(1/r);\\
 \partial_t u_k(r,0) & = \sum_{1\leq k_2 \leq \frac{\mu+\nu_k}{2}} B_{k,k_2} r^{-d-2\nu_k+2k_2} = r^{-\mu-\nu_k-1} Q_k(1/r).
\end{align*}
Here $P_k(z)$ and $Q_k(z)$ are polynomials as given in Proposition \ref{structure of PR}. Therefore we have 
\begin{align}
 \int_{\mathbb{S}^{d-1}} u_0 (r\theta) \Phi_k(\theta) d\theta & = r^{-\mu} P_k(1/r); \label{coefficient of u0}\\
 \int_{\mathbb{S}^{d-1}} u_1 (r\theta) \Phi_k(\theta) d\theta & = r^{-\mu-1} Q_k(1/r). \label{coefficient of u1}
\end{align}
Next we show the polynomials satisfy the inequalities in in Proposition \ref{structure of PR}. We have 
\begin{align*}
 \int_{\mathbb{S}^{d-1}} \nabla_\theta u_0 (r\theta) \nabla_\theta \Phi_k(\theta) d\theta & = - \int_{\mathbb{S}^{d-1}} u(r\theta,0) \Delta_{\mathbb{S}^{d-1}} \Phi_k(\theta) d\theta\\
 & = \nu_k (d-2+\nu_k)  r^{-\mu} P_k(1/r).
\end{align*}
Since $\nabla_\theta \Phi_k$ are orthogonal to each other with $L^2(\mathbb{S}^{d-1})$ norm $\sqrt{\nu_k(d-2+\nu_k)}$, we have 
\[
 \sum_{k=1}^\infty \nu_k (d-2+\nu_k) r^{-2\mu} \left|P_k(1/r)\right|^2 \leq  \|\nabla_\theta u_0(r\theta)\|_{L^2(\mathbb{S}^{d-1})}^2
\]
By the inequality $\|\nabla_\theta u_0 (r\theta)/r\|_{L^2\{x:|x|>R\}} = \|\slashed{\nabla} u_0\|_{L^2\{x:|x|>R\}} \leq \|\nabla u\|_{L^2\{x:|x|>R\}}$, we have 
\begin{equation} \label{upper bound P2}
\sum_{k=1}^\infty  \nu_k (d-2+\nu_k) \int_0^{1/R} \left|P_k(z)\right|^2 dz \lesssim_d \|\slashed{\nabla} u_0\|_{L^2(\{x: |x|>R\})}^2 < +\infty. 
\end{equation}
Similarly 
\[
 \sum_{k=0}^\infty \int_0^{1/R} \left|Q_k(z)\right|^2 dz = \|u_1\|_{L^2(\{x: |x|>R\})}^2 < +\infty.
\]
Next we differentiate \eqref{coefficient of u0} in $r$ and obtain 
\[
 \int_{\mathbb{S}^{d-1}} \partial_r u_0 (r\theta) \Phi_k(\theta) d\theta = r^{-\mu-1} (-\mu P_k(1/r) - (1/r) P'_k(1/r)).
\]
Following the same argument as above, we obtain 
\[
 \sum_{k=0}^\infty \int_0^{1/R} \left|\mu P_k(z) + z P'_k(z)\right|^2 dz \leq  \|\partial_r u_0\|_{L^2(\{x: |x|>R\})}^2 < +\infty. 
\]
Combining this inequality with \eqref{upper bound P2}, we have 
\[
 \sum_{k=0}^\infty \int_0^{1/R} \left( \nu_k (d-2+\nu_k) \left|P_k(z)\right|^2 + |zP'_k(z)|^2 \right) dz < +\infty.
\]
Since $\Phi_k(\theta)$ is a Hilbert basis, we may finally write $(u_0,u_1)$ in the following form by \eqref{coefficient of u0} and \eqref{coefficient of u1}. (These infinite sums are understood as convergence in $L^2(\mathbb{S}^{d-1})$ and $L^2(\{x: |x|>R\})$ respectively.)
\begin{align*}
 &u_0(r,\theta) = \sum_{k=0}^\infty r^{-\mu} P_k(1/r) \Phi_k(\theta)& &u_1(r,\theta) = \sum_{k=0}^\infty r^{-\mu-1} Q_k(1/r) \Phi_k(\theta).&
\end{align*}

\paragraph{Step two} Let us assume $P_k(z)$ and $Q_k(z)$ satisfy the conditions given in the proposition. Now we show $(u_0,u_1) \in P(R)$. We start by proving $(u_0,u_1) \in \dot{H}^1 \times L^2(\{x: |x|>R\})$. This is clear that the series 
\[
 \sum_{k=0}^\infty r^{-\mu-1} Q_k(1/r) \Phi_k(\theta)
\]
converges in the space $L^2(\{x:|x|>R\})$ and 
\begin{equation}
  \|u_1\|_{L^2(\{x: |x|>R\})}^2  =  \sum_{k=0}^\infty \int_0^{1/R} \left|Q_k(z)\right|^2 dz < +\infty.  \label{norm of u1}
\end{equation} 
Next we show $u_0 \in \dot{H}^1 (\{x:|x|>R\})$. We need the following technical lemma, whose proof is put in the Appendix. 
\begin{lemma} \label{polynomial lemma}
 Let $L\geq 2l >0$ and $P(z)$ be a polynomial of degree $\kappa$. Then we have 
\begin{align*}
 \max_{z \in [0,L]} |P(z)|^2 & \leq \frac{(\kappa+1)^2}{L} \int_0^L |P(z)|^2 dz;\\
 \int_0^l |zP'(z)|^2 dz & \leq \frac{2\kappa(\kappa+1) l}{L} \int_0^L |P(z)|^2 dz. 
\end{align*}  
\end{lemma}
\noindent As a result, we have 
\begin{align*}
 \left\|\sum_{k=N}^\infty r^{-\mu} P_k(1/r) \Phi_k(\theta)\right\|_{L^2(\mathbb{S}^{d-1})} & = r^{-\mu} \left(\sum_{k=N}^\infty |P_k(1/r)|^2\right)^{1/2} \\
 & \lesssim r^{-\mu} \left(\sum_{k=N}^\infty \frac{(\mu+\nu_k)^2 R}{2} \int_0^{1/R} |P_k(z)|^2 dz\right)^{1/2}.
\end{align*}
converges to zero uniformly in $r \in [R,+\infty)$ as $N\rightarrow +\infty$. Thus the series 
\[
 \sum_{k=0}^\infty r^{-\mu} P_k(1/r) \Phi_k(\theta) 
\]
converges to $u_0$ in $C([R,\infty); L^2(\mathbb{S}^{d-1}))$. Next we show 
\begin{align} 
  \nabla_x u_0 (r,\theta) & = \sum_{k=0}^\infty \nabla_x \left(r^{-\mu} P_k(1/r) \Phi_k(\theta)\right) \nonumber\\ 
  & = \sum_{k=0}^\infty r^{-\mu-1} \left\{P_k(1/r) \nabla_\theta \Phi_k(\theta) - [\mu P_k(1/r) + (1/r) P'_k(1/r)] \Phi_k(\theta) \theta\right\}. \label{full derivative}
\end{align}
Our assumption on $P_k(z)$, as well as the orthogonality of $\{\nabla_\theta \Phi_k\}_{k\geq 0}$ and $\{\Phi_k\}_{k\geq 0}$, guarantee that the series in the right hand side converges in $L^2([R,\infty) \times \mathbb{S}^{d-1}; r^{d-1} dr d\theta)$, or equivalently in $L^2(\{x: |x|>R\})$. Given any $\varphi \in C_0^\infty (\{x:|x|>R\})$, we have 
\[
\int_{|x|>R} \left( \sum_{k=0}^N r^{-\mu} P_k(1/r) \Phi_k(\theta) \right) \nabla_x \varphi (r, \theta ) dx  = - \int_{|x|>R}  \varphi(x)  \sum_{k=0}^N \nabla_x \left(r^{-\mu} P_k(1/r) \Phi_k(\theta)\right) dx.
\]
By the convergence of series we make $N \rightarrow +\infty$ and obtain
\[
 \int_{|x|>R} u_0(x) \nabla_x \varphi (x) dx
 = - \int_{|x|>R} \varphi(x) \sum_{k=0}^\infty \nabla_x \left(r^{-\mu} P_k(1/r) \Phi_k(\theta)\right) dx.
\]
This verifies \eqref{full derivative}. Please note that we always have $\nabla_\theta \Phi_k \cdot \theta = 0$, thus \eqref{full derivative} is actually an orthogonal decomposition. This immediately gives 
\begin{align}
 \|\partial_r u_0(r,\theta)\|_{L^2(\{x: |x|>R\})} & =   \sum_{k=0}^\infty \int_0^{1/R} \left|\mu P_k(z) + z P'_k(z)\right|^2 dz < +\infty; \label{norm of ur}\\
  \|\slashed{\nabla} u_0\|_{L^2(\{x: |x|>R\})}^2 & =  \sum_{k=1}^\infty \nu_k (d-2+\nu_k) \int_0^{1/R} \left|P_k(z)\right|^2 dz < +\infty. \label{norm of utheta}
\end{align}
In summary, we have $(u_0,u_1) \in \dot{H}^1 \times L^2(\{x: |x|>R\})$. For completeness we may define $u_0, u_1$ in the region $\{x: |x|\leq R\}$ so that $(u_0,u_1)\in \dot{H}^1 \times L^2(\Rm^d)$. Next we show that $\mathbf{S}_L (u_0,u_1)$ is a weakly non-radiative solution. First of all, if $1\leq k_1 \leq \frac{\mu+\nu_k}{2}$ we may find constants  $C_{k_1-1}, \cdots, C_1$ so that 
\[
 f(r,t) \doteq r^{-d-\nu_k+2k_1} + C_{k_1-1} t^2 r^{-d-\nu_k+2k_1-2} + \cdots + C_1 t^{2(k_1-1)} r^{-d-\nu_k+2}
\]
satisfies the equation $(\partial_t^2 - \partial_r^2 - \frac{d-1}{r}\partial_r) f(r,t) = -\frac{\nu_k(d+\nu_k-2)}{r^2} f(r,t)$. In fact these constant can be determined inductively. Therefore $v = f(r,t) \Phi_k(\theta)$ solves the equation 
\begin{equation}
 (\partial_t^2 - \Delta_x) v = (\partial_t^2 - \partial_r^2 - \frac{d-1}{r}\partial_r - \frac{\Delta_{\mathbf{S}^{d-1}}}{r^2}) v = 0, \qquad |x|=r>0. \label{cp2}
\end{equation}
with initial data $(r^{-d-\nu_k+2k_1}\Phi_k(\theta),0)$. A basic calculation shows that 
\[
 \lim_{t\rightarrow \pm \infty} \int_{|x|>|t|+R} |\nabla_{x,t} v(x,t)|^2 dx = 0.
\]
Similarly we may find a non-radiative solution $v$ to \eqref{cp2} with initial data $(0,r^{-d-\nu_k-1+2k_2} \Phi_k(\theta))$. By linearity, we may find a non-radiative solution $v_N$ to \eqref{cp2} with initial data $v_{0,N}, v_{1,N}$ so that 
\begin{align*}
 &v_{0,N} (r,\theta) = \sum_{k=0}^N r^{-\mu} P_k(1/r) \Phi_k(\theta);& &v_{1,N} (r,\theta) = \sum_{k=0}^N r^{-\mu-1} Q_k(1/r) \Phi_k(\theta).&
\end{align*}
By a standard centre cut-off and finite speed of propagation we obtain initial data $(u_{0,N}, u_{1,N})\in \dot{H}^1 \times L^2$ and corresponding free wave $u_N = \mathbf{S}_L (u_{0,N}, u_{1,N})$ so that 
\begin{align*}
 &v_{0,N} (r,\theta) = \sum_{k=0}^N r^{-\mu} P_k(1/r) \Phi_k(\theta);& &v_{1,N} (r,\theta) = \sum_{k=0}^N r^{-\mu-1} Q_k(1/r) \Phi_k(\theta);& &r>R&
\end{align*}
and 
\[
 \lim_{t\rightarrow \pm \infty} \int_{|x|>R+|t|} |\nabla_{x,t} u_N (x,t)|^2 dx = 0.
\]
By finite speed of energy propagation, we also have that $u = \mathbf{S}_L (u_0,u_1)$ satisfies 
\begin{align*}
 \limsup_{t\rightarrow \pm \infty} \int_{|x|>R+|t|} |\nabla_{x,t} (u (x,t)-u_N(x,t))|^2 dx \leq  \int_{|x|>R} \left(|\nabla u_0 - \nabla u_{0,N}|^2 + |u_1-u_{1,N}|^2\right) dx 
\end{align*}
We may combine the two limits above and obtain that the inequality 
\[
 \limsup_{t\rightarrow \pm \infty} \int_{|x|>R+|t|} |\nabla_{x,t} u (x,t)|^2 dx \lesssim_1 \int_{|x|>R} \left(|\nabla u_0 - \nabla u_{0,N}|^2 + |u_1-u_{1,N}|^2\right) dx. 
\]
holds for all $N \geq 1$. Finally we make $N\rightarrow +\infty$ and conclude that $(u_0,u_1)\in P(R)$.

\paragraph{Step three} Now we show that the identities and inequalities given in (i), (ii) and (iii) hold. Part (i) and (ii) have been proved in step two, see \eqref{norm of u1}, \eqref{full derivative}, \eqref{norm of ur} and \eqref{norm of utheta}. Now we consider part (iii). We have 
\[
 \partial_r u_0^\ast = \sum_{k=1}^\infty r^{-\mu-1} (-\mu P_k(1/r)-(1/r)P'_k(1/r)) \Phi_k(\theta)
\]
Thus 
\begin{align*}
 \int_{|x|>R_1} |\partial_r u_0^\ast|^2 dx = \sum_{k=1}^\infty \int_0^{1/R_1} |\mu P_k(z) + zP'_k(z)|^2 dz \lesssim_d \sum_{k=1}^\infty \int_0^{1/R_1} \left(|P_k(z)|^2 + |zP'_k(z)|^2\right) dz 
\end{align*}
We then apply Lemma \ref{polynomial lemma} and obtain
\begin{align*}
 \int_{|x|>R_1} |\partial_r u_0^\ast|^2 dx \lesssim_d \sum_{k=1}^\infty (\mu+\nu_k)^2 \frac{R}{R_1} \int_0^{1/R} |P_k(z)|^2 dz \lesssim_d (R/R_1) \int_{|x|>R} |\slashed{\nabla} u_0(x)|^2 dx.
\end{align*}
In order to find an upper bound of $\|\partial_r u_0\|_{L^2}$, we also need to consider the radial part $r^{-\mu} P_0(1/r) \Phi_0(\theta)$. In this case $\nu_0=0$ and $\Phi_0$ is simply a constant. We may follow the same argument above and obtain  
\begin{align*}
  \int_{|x|>R_1} \left|\partial_r [r^{-\mu} P_0(1/r) \Phi_0(\theta)]\right|^2 dx & =  \int_0^{1/R_1} |\mu P_0(z) + zP'_0(z)|^2 dz \\
  & \leq  \mu^2 \frac{R}{R_1} \int_0^{1/R} |\mu P_0(z) + zP'_0(z)|^2 dz\\
  & = \mu^2 \frac{R}{R_1} \int_{|x|>R} \left|\partial_r [r^{-\mu} P_0(1/r) \Phi_0(\theta)]\right|^2 dx\\
  & \lesssim_d \frac{R}{R_1} \int_{|x|>R} |\partial_r u_0(x)|^2 dx.
\end{align*}
Here we recall $\mu P_0(z) + zP'_0(z)$ is a polynomial of degree $\mu-1$ or less and apply Lemma \ref{polynomial lemma}. In summary, we use orthogonality to conclude ($R_1 \geq 2R$)
\begin{align*}
 \int_{|x|>R_1} |\partial_r u_0|^2 dx & = \int_{|x|>R_1} |\partial_r u_0^\ast|^2 dx + \int_{|x|>R_1} \left|\partial_r [r^{-\mu} P_0(1/r) \Phi_0(\theta)]\right|^2 dx\\
 & \lesssim \frac{R}{R_1} \int_{|x|>R} |\nabla u_0(x)|^2 dx. 
\end{align*}
This completes the proof of Proposition \ref{structure of PR}.
 
\subsection{Even dimensions} 

In this subsection we generalize our result on weakly non-radiative solutions to the even dimensions. 

\begin{proposition} \label{structure of PR even}
 Assume that $d\geq 2$ is an even integer and $\mu = d/2$. Then $(u_0,u_1) \in P(R)$ is equivalent to saying that there exist two sequences of polynomials $\{P_k(z)\}_{k\geq 0}$ and $\{Q_k(z)\}_{k\geq 0}$ of the following form 
 \begin{align*}
 &P_k(z) = \sum_{1\leq k_1 \leq \lfloor\frac{\mu+\nu_k}{2}\rfloor} A_{k,k_1} z^{\mu+\nu_k-2k_1};& &Q_k(z) = \sum_{1\leq k_2 \leq \lfloor\frac{\mu+\nu_k-1}{2}\rfloor} B_{k,k_2} z^{\mu+\nu_k-1-2k_2}.&
\end{align*}
with
\begin{align*}
 &\sum_{k=0}^\infty \int_0^{1/R} z\left( \nu_k (d-2+\nu_k) \left|P_k(z)\right|^2 + |zP'_k(z)|^2 \right) dz < +\infty;& &\sum_{k=0}^\infty \int_0^{1/R} z\left|Q_k(z)\right|^2 dz < +\infty.&
\end{align*}
so that 
 \begin{align}
 &u_0(r,\theta) = \sum_{k=0}^\infty r^{-\mu} P_k(1/r) \Phi_k(\theta)& &u_1(r,\theta) = \sum_{k=0}^\infty r^{-\mu-1} Q_k(1/r) \Phi_k(\theta).& \label{bh map}
\end{align}
Here the first identity holds for every $r>R$ in the sense of $L^2(\mathbb{S}^{d-1})$ convergence. the second one holds in the sense of $L^2(\{x: |x|>R\})$ convergence. In addition, we have 
\begin{itemize}
 \item[(i)] The derivative of $u_0$ can be given by 
 \begin{align*}
  \nabla_x u_0 (r,\theta) = \sum_{k=0}^\infty r^{-\mu-1} \left\{P_k(1/r) \nabla_\theta \Phi_k(\theta) - [\mu P_k(1/r) + (1/r) P'_k(1/r)] \Phi_k(\theta) \theta \right\}.
 \end{align*}
This identities holds in the sense of $L^2(\{x: |x|>R\})$ convergence. Here $\nabla_\theta \Phi_k(\theta)$ is in the tangent space of $\mathbb{S}^{d-1}$ at the point $\theta$ thus can be naturally embedded into $\Rm^d$. 
 \item[(ii)] The norms of $(u_0,u_1)$ can be determined by $P_k(z)$ and $Q_k(z)$'s:
 \begin{align*}
  \|\slashed{\nabla} u_0\|_{L^2(\{x: |x|>R\})}^2 & = \sum_{k=1}^\infty \nu_k (d-2+\nu_k) \int_0^{1/R} z\left|P_k(z)\right|^2 dz;\\
  \|u_1\|_{L^2(\{x: |x|>R\})}^2 & = \sum_{k=0}^\infty \int_0^{1/R} z \left|Q_k(z)\right|^2 dz; \\
  \|\partial_r u_0\|_{L^2(\{x: |x|>R\})} & = \sum_{k=0}^\infty \int_0^{1/R} z |zP'_k(z)+\mu P_k(z)|^2 dz < +\infty.
 \end{align*}. 
 \item[(iii)] The derivative $\partial_r u_0$ satisfies the following decay estimates ($R_1 \geq 2R$)
 \begin{align*}
  \int_{|x|>R_1} |\partial_r u_0(x)|^2 dx & \lesssim (R/R_1) \int_{|x|>R} |\nabla u_0(x)|^2  dx;\\
  \int_{|x|>R_1} |\partial_r u_0^\ast (x)|^2 dx  & \lesssim (R/R_1) \int_{|x|>R} |\slashed{\nabla} u_0(x)|^2 dx. 
 \end{align*}
 Here $u_0^\ast$ is the non-radial part of $u_0$ defined by $u_0^\ast = u_0 - r^{-\mu} P_0(1/r) \Phi_0$. 
\end{itemize}
\end{proposition}
\noindent The proof in the even dimensions is almost the same as in the odd dimensions thus we omit it here. The main difference is that we rely on a slightly modified version of the technical lemma about polynomials, which is given below and proved in the appendix. 
\begin{lemma} \label{polynomial lemma even}
 Let $L\geq 2l >0$ and $P(z)$ be a polynomial of degree $\kappa$. Then we have 
\begin{align*}
 \max_{z \in [0,L]} z |P(z)|^2 & \leq \frac{2(\kappa+1)^2}{L} \int_0^L z |P(z)|^2 dz;\\
 \int_0^l z |zP'(z)|^2 dz & \leq \frac{2\kappa(\kappa+2) l}{L} \int_0^L z |P(z)|^2 dz. 
\end{align*}  
\end{lemma}

\section{Non-linear Non-radiative Solutions} \label{sec: nonlinear case}

In this section we show that non-radiative solutions to a wide range of nonlinear wave equations in the three-dimensional case share the same asymptotic behaviour as non-radiative free waves, without the radial assumption. 

\paragraph{Assumptions} We consider the energy-critical non-linear wave equation in $\Rm^3$
\[
 \partial_t^2 u - \Delta u = F(x,t,u), \qquad (x,t) \in \Rm^3 \times \Rm.
\]
Here the nonlinear term $F(x,t,u)$ satisfies 
\begin{align}
 &|F(x,t,u)| \leq C |u|^5; & &|F(x,t,u_1) - F(x,t,u_2)| \leq C(|u_1|^4 + |u_2|^4)|u_1-u_2|.& \label{nonlinear assumption}
\end{align} 
This covers both the defocusing ($F(x,t,u) = -|u|^4 u$) and focusing ($F(x,t,u) = |u|^4 u$) wave equations, which have been extensively studied in the past decades. 

\subsection{Preliminary results}

We first give a few preliminary results and introduce a few notations. 

\paragraph{Radiation fields} Radiation field describes the asymptotic behaviour of free waves as time tends to infinity. In its earlier history radiation field was mainly a conception in mathematical physics. See Friedlander \cite{radiation1, radiation2}, for instance. The following modern version is given in \cite{dkm3}. 

\begin{theorem}[Radiation fields] \label{radiation}
Assume that $d\geq 3$ and let $u$ be a solution to the free wave equation $\partial_t^2 u - \Delta u = 0$ with initial data $(u_0,u_1) \in \dot{H}^1 \times L^2(\Rm^d)$. Then ($u_r$ is the derivative in the radial direction)
\[
 \lim_{t\rightarrow \pm \infty} \int_{\Rm^d} \left(|\nabla u(x,t)|^2 - |u_r(x,t)|^2 + \frac{|u(x,t)|^2}{|x|^2}\right) dx = 0
\]
 and there exist two functions $G_\pm \in L^2(\Rm \times \mathbb{S}^{d-1})$ so that
\begin{align*}
 \lim_{t\rightarrow \pm\infty} \int_0^\infty \int_{\mathbb{S}^{d-1}} \left|r^{\frac{d-1}{2}} \partial_t u(r\theta, t) - G_\pm (r\mp t, \theta)\right|^2 d\theta dr &= 0;\\
 \lim_{t\rightarrow \pm\infty} \int_0^\infty \int_{\mathbb{S}^{d-1}} \left|r^{\frac{d-1}{2}} \partial_r u(r\theta, t) \pm G_\pm (r\mp t, \theta)\right|^2 d\theta dr & = 0.
\end{align*}
In addition, the maps $(u_0,u_1) \rightarrow \sqrt{2} G_\pm$ are bijective isometries from $\dot{H}^1 \times L^2(\Rm^d)$ to $L^2 (\Rm \times \mathbb{S}^{d-1})$. 
\end{theorem}

\noindent We call $G_\pm$ radiation fields associated to the free wave $u$. Throughout this section we utilize the notations $\mathbf{T}_\pm$ for the linear map from the initial data $(u_0,u_1)$ to the corresponding radiation fields $G_\pm$. It immediately follows the theorem that 
\[
 \lim_{t\rightarrow \pm \infty} \int_{|x|>R+|t|} |\nabla_{t,x} u(x,t)|^2 dx = 2 \int_R^\infty \int_{\mathbb{S}^2} |G_\pm (s, \theta)|^2 d\theta ds. 
\]
In addition, the map between $G_\pm$ is an isometry given explicitly by 
\[
 G_+(s,\theta) = \left\{\begin{array}{ll} (-1)^{\frac{d-1}{2}} G_-(-s,-\theta), & d\; \hbox{is odd;}\\ (-1)^{\frac{d}{2}} (\mathcal{H} G_-)(-s,-\theta), & d\; \hbox{is even}.\end{array}\right.
\]
This can proved in different methods. Please refer to C\^{o}te-Laurent \cite{newradiation}, Duyckaerts-Kenig-Merle \cite{dkmnonradial} and Li-Shen-Wei \cite{shenradiation}, for examples. As a result, the following identity holds for all odd dimensions $d\geq 3$:
\begin{equation} \label{exterior energy identity}
 \sum_{\pm} \lim_{t\rightarrow \pm \infty} \int_{|x|>R+|t|} |\nabla_{t,x} u(x,t)|^2 dx = 2 \int_{|r|>R} \int_{\mathbb{S}^2} |G_- (s, \theta)|^2 d\theta ds.
\end{equation}
As a result, the $L^2$ decay rate of radiation field $G_- (s,\theta)$ near the infinity indicates to what extent the free wave $u$ looks like a non-radiative solution. 

\paragraph{Decay of linear non-radiative solutions} Another important ingredient of our estimate on non-linear non-radiative solutions is the corresponding decay estimates of linear non-radiative solutions. We claim that given any constant $\kappa \in (0,1/5)$, the following inequality holds 
\begin{equation} \label{decay estimate assumption}
 \|u\|_{L_t^5 L^{10}(\{x: |x|>r+|t|\})} \lesssim_\kappa (R/r)^\kappa E^{1/2},
\end{equation}
for any $r \geq R > 0$ and $R$-weakly non-radiative linear wave $u$ with a finite energy $E$, i.e. a finite-energy solution to the homogeneous linear wave equation $\partial_t^2 u - \Delta u = 0$ so that 
\[
 \lim_{t\rightarrow \pm \infty} \int_{|x|>R+|t|} |\nabla_{t,x} u(x,t)|^2 dx = 0.
\]
In fact, it was prove in Li-Shen-Wang \cite{3Dnonradiativedecay} that any $R$-weakly non-radiative linear wave $u$ satisfies that inequality 
\begin{equation} \label{original decay estimate}
  \|u\|_{L_t^\infty L^{6}(\{x: |x|>r+|t|\})} \lesssim (R/r)^{1/3} E^{1/2}. 
\end{equation}
We may interpolate it with a regular Strichartz estimate (see Ginibre-Velo \cite{strichartz})
\[
 \|u\|_{L_t^p L_x^q (\Rm \times \Rm^3)} \lesssim_{p,q} E^{1/2}
\]
with $p = 2^+$ and $q= \infty^-$ and conclude that the inequality \eqref{decay estimate assumption} holds for any $\kappa \in (0,1/5)$. 

\subsection{Statement and Proof}

\begin{proposition}
Let $u$ be an $R$-weakly non-radiative solution to the non-linear wave equation 
\[
 \left\{\begin{array}{ll} \partial_t^2 u - \Delta u = F(x,t,u), & (x,t)\in \Rm^3 \times \Rm; \\ (u,u_t)|_{t=0} = (u_0,u_1) \in \dot{H}^1\times L^2(\Rm^3). & \end{array} \right.
\]
Here the nonlinear term satisfies \eqref{nonlinear assumption}. Then we have
\begin{itemize}
 \item[(a)] Given any $\kappa \in (0,1/5)$, the radiation field $G_-(s,\omega)$ associated to the linear wave $\mathbf{S}_L (u_0,u_1)$ satisfies a decay estimate 
 \[
  \|G_-\|_{L^2 (\{s: |s|>r\}\times \mathbb{S}^2)} \lesssim r^{-5\kappa}, \qquad \forall r\gg R.
 \]
 It is equivalent to saying (see \eqref{exterior energy identity})
 \[
  \lim_{t\rightarrow \pm \infty} \int_{|x|>r+|t|} |\nabla_{t,x} \mathbf{S}_L(u_0,u_1) (x,t)|^2 dx \lesssim r^{-10\kappa}, \qquad \forall r\gg R.
 \]
 \item[(b)] The initial data $u_0$ satisfy the decay estimate
 \[
  \int_{|x|>r} |\partial_r u_0(x)|^2 dx \lesssim r^{-1}, \qquad \forall r\gg R.
 \]
 \item[(c)] We also the decay estimate
 \[
  \sup_{t\in \Rm} \int_{|x|>r+|t|} |u(x,t)|^6 dx \lesssim r^{-2}, \qquad \forall r\gg R.
 \]
\end{itemize}
\end{proposition}
\begin{proof}
 Let us first introduce a notation for convenience. We define
\[
 S(r) = \|G_-\|_{L^2 (\{s: |s|>r\}\times \mathbb{S}^2)} = \left(\int_{|s|>r} \int_{\mathbb{S}^2} |G_-(s,\omega)|^2 d\omega ds\right)^{1/2}. 
\]
Given any $r \gg r_1 \gg R$, we may break $G_-$ into two parts 
\begin{align*}
 &G_1 (s,\omega) = \left\{\begin{array}{ll} G_-(s,\omega), & |s|\leq r_1; \\ 0, & |s|>r_1; \end{array} \right. & &G_2 (s,\omega) = \left\{\begin{array}{ll} 0, & |s|\leq r_1; \\ G_-(s,\omega), & |s|>r_1. \end{array} \right.&
\end{align*}
Therefore we have 
\begin{equation}
 (u_0,u_1) = \mathbf{T}_-^{-1} G_1 + \mathbf{T}_-^{-1} G_2.
\end{equation} 
We also define $\chi_r (x,t)$ to be the characteristic function of the exterior region $\Omega(r) = \{(x,t): |x|>|t|+r\}$ and 
\[
 \|v\|_{Y(r)} = \|\chi_r (x,t) u\|_{L^5 L^{10} (\Rm \times \Rm^3)} = \|v\|_{L_t^5 L^{10}(\{x: |x|>r +|t|\})},
\]
Next we give a reasonable upper bound of $\|\mathbf{S}_L(u_0,u_1)\|_{Y(r)}$ by our decay estimate assumption. In fact we have 
\begin{align}
 \|\mathbf{S}_L (u_0,u_1)\|_{Y(r)} & \leq \|\mathbf{S}_L \mathbf{T}_-^{-1} G_1\|_{Y(r)} + \|\mathbf{S}_L \mathbf{T}_-^{-1} G_2\|_{Y(r)} \nonumber \\
 & \lesssim (r_1/r)^\kappa \|G_1\|_{L^2} + \|G_2\|_{L^2}\nonumber \\
 & \lesssim (r_1/r)^{\kappa}  + S(r_1). \label{upper bound on Y part 1}
\end{align}
Here we utilize the fact that $G_1$ is supported in $[-r_1,r_1]\times \mathbb{S}^2$ thus the linear free wave $\mathbf{S}_L \mathbf{T}_-^{-1} G_1$ with radiation field $G_1$ is an $r_1$-weakly non-radiative free wave. We then apply \eqref{decay estimate assumption} on the $G_1$ part and the classic Strichartz estimate on the $G_2$ part. Now we consider a modified non-linear wave equation 
\begin{equation} \label{modified NLW} 
 \left\{\begin{array}{ll} \partial_t^2 v - \Delta v = \chi_r (x,t) F(x,t,v), & (x,t)\in \Rm^3 \times \Rm; \\ (v,v_t)|_{t=0} = (u_0,u_1) \in \dot{H}^1\times L^2(\Rm^3). & \end{array} \right.
\end{equation}
First of all, the following inequalities hold by our assumption on the nonlinear term $F$. 
\begin{align*}
 \left\|\chi_r F(x,t,v)\right\|_{L^1 L^2 (\Rm \times \Rm^3)} & \lesssim \|v\|_{Y(r)}^5; \\
 \left\|\chi_r F(x,t,v_1) - \chi_r F(x,t,v_2)\right\|_{L^1 L^2(\Rm \times \Rm^3)} & \lesssim (\|v_1\|_{Y(r)}^4 + \|v_2\|_{Y(r)}^4)\|v_1 - v_2\|_{Y(r)}.
\end{align*}
We also recall the classic Strichartz estimate (see \cite{strichartz}): if $w$ solves the 3D linear wave equation $\partial_t^2 w - \Delta w = F$ with initial data $(w_0,w_1)$, then 
\[
 \|w\|_{L^5 L^{10} (\Rm \times \Rm^3)} + \|(w,w_t)\|_{C(\Rm_t; \dot{H}^1 \times L^2)}  \lesssim \|(w_0,w_1)\|_{\dot{H}^1 \times L^2} + \|F\|_{L^1 L^2(\Rm \times \Rm^3)}. 
\]
We may combine all these inequalities, apply a standard fixed-point argument of contraction map and conclude that as long as $\|\mathbf{S}_L (u_0,u_1)\|_{Y(r)}$ is sufficiently small, which holds under our assumption $r \gg r_1 \gg R$ by \eqref{upper bound on Y part 1}, the equation $\eqref{modified NLW}$ always has a global-in-time solution $v$, so that 
\begin{equation} \label{small data theory bound}
 \|v\|_{Y(r)} \leq 2 \|\mathbf{S}_L (u_0,u_1)\|_{Y(r)}.
\end{equation}
More details about the fixed-point argument of this kind can be found, for instance, in Pecher \cite{Pecher}. Furthermore, we may write $v$ as a sum of two terms
\[
 v= v_1 + v_2.
\]
They are the linear propagation part and the contribution of non-linear term, respectively:
\begin{align}
 &v_1 = \mathbf{S}_L (u_0,u_1);& & v_2 =  \int_0^t \frac{\sin (t-\tau)\sqrt{-\Delta}}{\sqrt{-\Delta}} (\chi_r F(\cdot, \tau, v(\cdot, \tau))) d\tau.& \label{decomposition of v}
\end{align}
The triangle inequality in $L^2$ space gives 
\begin{align*}
 \left(\int_{|x|>r +|t|} |\nabla_{t,x} v|^2 dx\right)^{1/2} \geq  \left(\int_{|x|>r +|t|} |\nabla_{t,x} v_1|^2 dx\right)^{1/2} - \left(\int_{|x|> r+|t|} |\nabla_{t,x} v_2|^2 dx\right)^{1/2}
\end{align*}
for any given time $t$.  A comparison of our modified non-linear wave equation $\eqref{modified NLW}$ with the original one shows that $u(x,t) \equiv v(x,t)$ in the exterior region $\Omega(r)$ by finite speed of propagation. Therefore our non-radiative assumption on $u$ also applies on $v$ in the exterior region $\Omega(r)$. This gives 
\[
 \lim_{t\rightarrow \pm \infty} \int_{|x|>r +|t|} |\nabla_{t,x} v|^2 dx = 0.
\]
Therefore we have 
\[
 \liminf_{t\rightarrow \pm \infty} \int_{|x|>r+|t|} |\nabla_{t,x} v_2(x,t)|^2 dx \geq \lim_{t\rightarrow \pm \infty} \int_{|x|>r+|t|} |\nabla_{t,x} v_1(x,t)|^2 dx.
\]
We then recall the property of radiation field and obtain 
\[
\sum_{\pm} \lim_{t\rightarrow \pm \infty} \int_{|x|>r+|t|} |\nabla_{t,x} v_1(x,t)|^2 dx = 2 \int_{|s|>r} \int_{\mathbb{S}^2} |G_-(s,\omega)|^2 d\omega ds = 2 S^2 (r).
\]
We may also find an upper bound of the integral about $v_2$ by Strichartz estimates
\begin{align*}
 \int_{|x|>r+|t|} |\nabla_{t,x} v_2(x,t)|^2 dx & \leq \int_{\Rm^3} |\nabla_{t,x} v_2(x,t)|^2 dx \\
 & \leq \left\|\chi_r  F(x,t,v)\right\|_{L^1 L^2 (\Rm \times \Rm^3)}^2\\
 & \lesssim \|v\|_{Y(r)}^{10}.
\end{align*}
Combining these inequalities we obtain $S(r)\lesssim \|v\|_{Y(r)}^5$. We then utilize the upper bound given in \eqref{small data theory bound} and obtain   
\begin{equation} \label{S by SL}
 S(r) \lesssim \|\mathbf{S}_L (u_0,u_1)\|_{Y(r)}^5, \qquad r\gg R. 
\end{equation} 
A combination of this inequality with \eqref{upper bound on Y part 1} immediately gives a recursion formula when $r \gg r_1 \gg R$.
\[
 S(r) \lesssim (r_1/r)^{5\kappa} + S^5(r_1).
\]
We then apply Lemma \ref{recursion lemma}, whose statement and proof is postponed to the appendix, and conclude that given any $\beta \in (0, 5\kappa)$, the following estimate holds if $r \geq R_0(u,\kappa,\beta)$ is sufficiently large
\[
 S(r) \leq r^{-\beta}.
\]
Next we give a more detailed estimate of $\|\mathbf{S}_L (u_0,u_1)\|_{Y(r)}$ as $r\rightarrow +\infty$. We fix a constant $\beta \in (\kappa, 5\kappa)$, choose $R_0 = R_0(u, \kappa, \beta)$ accordingly as above and define 
\begin{align*}
 &G_0(s,\omega) = \left\{\begin{array}{ll} G_-(s,\omega), & |s|\leq R_0; \\ 0, & |s|>R_0; \end{array} \right.& &G_j(s,\omega) = \left\{\begin{array}{ll} G_-(s,\omega), & 2^{j-1} R_0 < |s|\leq 2^j R_0; \\ 0, & \hbox{otherwise}; \end{array} \right.\quad j\geq 1.&
\end{align*} 
Thus we have 
\[
 (u_0,u_1) = \sum_{j=0}^\infty \mathbf{T}_-^{-1} G_j.
\]
If $r \in [2^n R_0, 2^{n+1} R_0]$ for an integer $n \geq 0$, then we have 
\begin{align*}
 \|\mathbf{S}_L(u_0,u_1)\|_{Y(r)} & \leq \sum_{j=0}^n \|\mathbf{S}_L \mathbf{T}_-^{-1} G_j\|_{Y(r)} + \left\|\mathbf{S}_L \mathbf{T}_-^{-1} \left(\sum_{j=n+1}^\infty G_j \right)\right\|_{Y(r)}\\
 & \lesssim \sum_{j=0}^n (2^j R_0/r)^{\kappa} \|G_j\|_{L^2} + \left\|\sum_{j=n+1}^\infty G_j\right\|_{L^2}\\
 & \lesssim r^{-\kappa} + \sum_{j=1}^{n} 2^{\kappa j} r^{-\kappa} S(2^{j-1} R_0) + S(2^n R_0)\\
 & \lesssim r^{-\kappa} + \sum_{j=1}^{n} 2^{\kappa j} r^{-\kappa} (2^{j-1} R_0)^{-\beta} + (2^n R_0)^{-\beta}\\
 & \lesssim r^{-\kappa}.
\end{align*}
We apply the decay estimate \eqref{decay estimate assumption} and use the upper bound $S(r) \leq r^{-\beta}$ here. Finally we recall \eqref{S by SL} and conclude that the inequality $S(r) \lesssim r^{-5\kappa}$ holds if $r \geq R_1$ is sufficiently large. This finishes the proof of part (a). The proof of part (b) is similar the final stage of proof for part (a). We first fix a constant $\kappa \in (1/10,1/5)$. According to part (a), there exists $R_1>R$ so that $S(r) \lesssim r^{-5\kappa}$ holds for $r\geq R_1$. We define 
\begin{align*}
 &G_0(s,\omega) = \left\{\begin{array}{ll} G_-(s,\omega), & |s|\leq R_1; \\ 0, & |s|>R_1; \end{array} \right.& &G_j(s,\omega) = \left\{\begin{array}{ll} G_-(s,\omega), & 2^{j-1} R_1 < |s|\leq 2^j R_1; \\ 0, & \hbox{otherwise}; \end{array} \right.\quad j\geq 1;&
\end{align*} 
and 
\begin{align*}
 &(u_0, u_1) = \sum_{j=0}^\infty (u_{0,j}, u_{1,j}),& &(u_{0,j}, u_{1,j}) = \mathbf{T}_-^{-1} G_j.&
\end{align*}
Since $(u_{0,j}, u_{1,j}) \in P(2^j R_1)$, if $r > 2^j R_1$, then we may apply Proposition \ref{structure of PR} and obtain 
\begin{align*} 
 \left(\int_{|x|>r} |\partial_r u_{0,j}(x)|^2 dx \right)^{1/2} & \lesssim (2^j R_1 /r)^{1/2} \left(\int_{|x|>2^j R_1} |\nabla u_{0,j}(x)|^2 dx\right)^{1/2} \\
 & \lesssim (2^j R_1 /r)^{1/2} \|G_j\|_{L^2}.
\end{align*}
Furthermore, if we also have $j\geq 1$, then we may use the upper bound of $S(r)$ and obtain $\|G_j\|_{L^2} \leq S(2^{j-1} R_1) \lesssim (2^{j-1} R_1)^{-5\kappa}$. As a result, we have
\[
 \left(\int_{|x|>r} |\partial_r u_{0,j}(x)|^2 dx \right)^{1/2} \lesssim (2^j R_1 /r)^{1/2} (2^{j-1} R_1)^{-5\kappa} \lesssim (2^j R_1)^{1/2-5\kappa} r^{-1/2}.
\]
Now we assume $r > R_1$. Thus there exists $n \geq 1$ so that $2^{n-1} R_1 < r \leq 2^n R_1$. By the upper bounds given above, we have
\begin{align*}
 \left(\int_{|x|>r} |\partial_r u_0 (x)|^2 dx \right)^{1/2} & \leq  \left(\int_{|x|>r} |\partial_r u_{0,0} (x)|^2 dx \right)^{1/2} + \sum_{j=1}^{n-1}  \left(\int_{|x|>r} |\partial_r u_{0,j} (x)|^2 dx \right)^{1/2} \\
 & \qquad +  \left(\int_{|x|>r} \left|\partial_r \sum_{j=n}^\infty u_{0,j} (x)\right|^2 dx \right)^{1/2}\\
 & \lesssim (R_1/r)^{1/2} \|G_0\|_{L^2} + \sum_{j=1}^{n-1} (2^j R_1)^{1/2-5\kappa} r^{-1/2} + \left\|\sum_{j=n}^\infty G_j\right\|_{L^2}\\
 & \lesssim r^{-1/2} + S(2^{n-1} R_1)\\
 & \lesssim r^{-1/2} + (2^{n-1} R_1)^{-5\kappa}\\
 & \lesssim r^{-1/2}.
\end{align*}
This finishes the proof of part (b). Finally we prove part (c). We fix $\kappa \in (1/10,1/5)$ and use the same decomposition of $G_-$ and $(u_0,u_1)$ as in the proof of part (b). We apply the decay estimate \eqref{original decay estimate} and obtain that if $r>2^j R_1$, then
\begin{align*}
 \sup_{t\in \Rm} \left(\int_{|x|>r+|t|}|\mathbf{S}_L (u_{0,j}, u_{1,j})|^6 dx \right)^{1/6} \lesssim (2^j R_1/r)^{1/3} \|G_j\|_{L^2}.
\end{align*}
As a result, if $2^{n-1} R_1 < r < 2^n R_1$, then we have 
\begin{align*}
 \sup_{t\in \Rm} \left(\int_{|x|>r+|t|}|\mathbf{S}_L (u_{0}, u_{1})|^6 dx \right)^{1/6} & \leq \sum_{j=0}^{n-1} \sup_{t\in \Rm} \left(\int_{|x|>r+|t|}|\mathbf{S}_L (u_{0,j}, u_{1,j})|^6 dx \right)^{1/6} \\
 & \qquad + \sup_{t\in \Rm} \left(\int_{|x|>r+|t|}\left|\mathbf{S}_L \left(\sum_{j=n}^\infty (u_{0,j}, u_{1,j})\right) (x,t) \right|^6 dx \right)^{1/6}\\
 & \lesssim \sum_{j=0}^{n-1} (2^j R_1/r)^{1/3} \|G_j\|_{L^2} + \left\|\sum_{j=n}^\infty G_j\right\|_{L^2}.
\end{align*}
We then apply the $L^2$ decay estimate of $G_-$ given in part (a) and obtain ($r>R_1$)
\begin{align*}
 \sup_{t\in \Rm} \left(\int_{|x|>r+|t|}|\mathbf{S}_L (u_{0}, u_{1})|^6 dx \right)^{1/6} & \lesssim (R_1/r)^{1/3} + \sum_{j=1}^{n-1} (2^j R_1/r)^{1/3} (2^{j-1} R_1)^{-5\kappa} + (2^{n-1} R_1)^{-5\kappa} \\
 & \lesssim r^{-1/3} + \sum_{j=1}^{n-1} (2^j R_1)^{1/3 - 5\kappa} r^{-1/3} + r^{-5\kappa}\\
 & \lesssim r^{-1/3}.  
\end{align*}
Next we recall that if we let $v$ solves \eqref{modified NLW} and define $v_1, v_2$ accordingly as in \eqref{decomposition of v}, then 
\[
 u(x,t) = v(x,t) = v_1(x,t) + v_2(x,t)
\]
holds in the exterior region $\{(x,t): |x|>r+|t|\}$. Our argument above has already given $L^6$ upper bound of $v_1 = \mathbf{S}_L (u_0,u_1)$. It suffices to consider the upper bound of $v_2$. By the Strichartz estimates, we have
\begin{align*}
 \sup_{t\in \Rm} \|v_2(\cdot, t)\|_{L^6 (\Rm^3)} \lesssim \sup_{t\in \Rm} \|v_2(\cdot,t)\|_{\dot{H}^1(\Rm^3)} \lesssim \|\chi_r F(x,t,v)\|_{L^1 L^2(\Rm \times \Rm^3)} \lesssim \|v\|_{Y(r)}^5
\end{align*}
Finally we recall \eqref{small data theory bound} and the estimate $\|\mathbf{S}_L (u_0,u_1)\|_{Y(r)} \lesssim r^{-\kappa}$ given in part (a), if $r$ is sufficiently large, and obtain $\|v\|_{Y(r)} \lesssim r^{-\kappa}$. Combining this with the inequality above we have 
\[
 \sup_{t\in \Rm} \|v_2(\cdot, t)\|_{L^6 (\Rm^3)} \lesssim r^{-5\kappa}, \qquad r\gg R. 
\]
We collect upper bounds of $v_1 = \mathbf{S}_L(u_0,u_1)$ and $v_2$ to conclude the proof of part (c).
\[
 \sup_{t\in \Rm} \left(\int_{|x|>r+|t|}|u(x,t)|^6 dx \right)^{1/6} \lesssim r^{-1/3}, \qquad r\gg R. 
\]
\end{proof}

\section{Appendix}
In this section we prove a few technical lemmata. The authors believe that these results are probably previously known. For completeness we still give their proof. 
\paragraph{Polynomial estimates} We start by Lemma \ref{polynomial lemma}. By change of variables $x = 2z/L-1$, we may rewrite this technical lemma as below. 
\begin{lemma} \label{polynomial lemma 2}
Let $0<\delta \leq 1$ and $P(x)$ be a polynomial of degree $\kappa$. Then we have 
\begin{align*}
 \max_{x \in [-1,1]} |P(x)|^2 & \leq \frac{(\kappa+1)^2}{2} \int_{-1}^1 |P(x)|^2 dx;\\
 \int_{-1}^{-1+\delta} |(x+1)P'(x)|^2 dx & \leq \kappa(\kappa+1) \delta \int_{-1}^1 |P(x)|^2 dx. 
\end{align*}  
\end{lemma}
\begin{proof}
 Let us recall Legendre polynomials $P_n$ defined by
\[
 P_n(x) = \frac{1}{2^n n!} \frac{d^n}{dx^n} (x^2-1)^n.
\]
It is well known that $\{P_n\}_{n=0,1,2,\cdots}$ are orthogonal to each other in $L^2([-1,+1])$ with norm $\|P_n\|_{L^2}^2 = \frac{2}{2n+1}$. In addition, these polynomials satisfy $|P_n(x)|\leq 1, \forall |x|\leq 1$ and the differential equation
\[
 \frac{d}{dx}\left[(1-x^2)\frac{d}{dx}P_n(x)\right] + n(n+1) P_n(x) = 0. 
\]
More details about the properties of Legendre polynomials can be found, for instance, in Folland \cite{fourierappli}. We consider the orthogonal decomposition of $P(x)$:
\[
  P(x) = \sum_{n=0}^\kappa a_n P_n(x) \qquad \Rightarrow \qquad \int_{-1}^1 |P(x)|^2 dx = \sum_{n=0}^\kappa \frac{2|a_n|^2}{2n+1}.
\]
This immediately gives 
\begin{align*}
 \max_{x\in [-1,1]} |P(x)|^2 \leq \left(\sum_{n=0}^\kappa |a_n|\right)^2 \leq \left(\sum_{n=0}^{\kappa} \frac{2n+1}{2}\right)\left(\sum_{n=0}^\kappa \frac{2|a_n|^2}{2n+1}\right)= \frac{(\kappa+1)^2}{2}\int_{-1}^1 |P(x)|^2 dx 
\end{align*}
We also have 
\[
 \int_{-1}^{-1+\delta} |(x+1)P'(x)|^2 dx  \leq \delta \int_{-1}^{-1+\delta} (1-x^2) |P'(x)|^2 dx \leq \delta \int_{-1}^{1} (1-x^2) |P'(x)|^2 dx 
\]
We then integrate by parts, use the differential equation above and obtain 
\begin{align*}
 \int_{-1}^{1} (1-x^2) |P'(x)|^2 dx & = - \int_{-1}^1 P(x) \cdot \frac{d}{dx} [(1-x^2) P'(x)] dx \\
 & = \int_{-1}^1 \left(\sum_{n=0}^\kappa a_n P_n(x)\right) \left(\sum_{k=0}^\kappa n(n+1) a_n P_n(x) \right) dx\\
 & = \sum_{n=0}^\kappa \frac{2n(n+1)|a_n|^2}{2n+1}\\
 &  \leq \kappa(\kappa+1) \int_{-1}^1 |P(x)|^2 dx.
\end{align*}
Combining these two inequalities, we finish the proof. 
\end{proof}
\noindent We also need a similar lemma, where $dx$ is substituted by $(x+1) dx$. This immediately gives Lemma \ref{polynomial lemma even} by a change of variables $x = 2z/L-1$.
\begin{lemma} \label{polynomial lemma even 2}
Let $0<\delta \leq 1$ and $P(x)$ be a polynomial of degree $\kappa$. Then we have 
\begin{align}
 \max_{x \in [-1,1]} (x+1) |P(x)|^2 & \leq (\kappa+1)^2 \int_{-1}^1 (x+1) |P(x)|^2 dx; \label{first inequality} \\
 \int_{-1}^{-1+\delta} (x+1)^3 |P'(x)|^2 dx & \leq \kappa(\kappa+2) \delta \int_{-1}^1 (x+1) |P(x)|^2 dx. \label{second inequality}
\end{align}  
\end{lemma}
\begin{proof}
We define $Q_n(x)$ to be the modified Legendre polynomial of degree $n$:
\[ 
 Q_n(x) = \frac{1}{2^{n+1} (n+1)!}\frac{d^{n+1}}{dx^{n+1}} [(x+1)^n (x-1)^{n+1}] = \frac{(2n+1)!}{2^{n+1} n! (n+1)!} x^{n} + \cdots.
\]
If $n \geq m$ are nonnegative integers, then we may apply integration by parts and obtain
\begin{align*}
 \int_{-1}^1 (x+1) Q_n (x) Q_m(x) dx & = \frac{(-1)^{n+1}}{2^{n+1} (n+1)!} \int_{-1}^1 (x+1)^n (x-1)^{n+1} \frac{d^{n+1}}{dx^{n+1}} \left[(x+1) Q_m(x)\right] dx.
\end{align*}
A basic calculation shows 
\[
 \frac{d^{n+1}}{dx^{n+1}} \left[(x+1) Q_m(x)\right] = \left\{\begin{array}{ll} \frac{(2n+1)!}{2^{n+1} n!}, & \hbox{if} \; m=n;\\ 0, & \hbox{if}\; m<n. \end{array}\right.
\]
Therefore $\{Q_n (x)\}_{n\geq 0}$ are orthogonal to each other in the Hilbert space $L^2 ([-1,1]; (x+1)dx)$ and the norms of these polynomials are given by 
\[
 \|Q_n\|_{L^2([-1,1]; (x+1)dx)}^2 = \frac{1}{2(n+1)}. 
\]
In addition, these polynomials satisfy a similar differential equation to Legendre polynomials. 
\begin{equation} \label{diff equation mod L}
 \frac{d}{dx} \left[(x+1)(1-x^2) \frac{d}{dx} Q_n (x) \right] + n(n+2)(x+1) Q_n(x) = 0.
\end{equation}
In order to prove this identity, we observe that $\frac{d}{dx}\left[(x+1)(x^2-1) \frac{d}{dx} Q_n (x) \right]$ is a polynomial of degree $n+1$ and contain a factor of $x+1$. Thus we may write 
\[
 \frac{d}{dx}\left[(x+1)(x^2-1) \frac{d}{dx} Q_n (x) \right] = \sum_{j=0}^n a_j (x+1) Q_j (x).
\]
We multiply both sides by $Q_j(x)$, integrate from $x=-1$ to $x=1$ and apply integration by parts 
\begin{align*}
 \frac{a_j}{2(j+1)} & = \int_{-1}^1 Q_j(x) \frac{d}{dx}\left[(x+1)(x^2-1) \frac{d}{dx} Q_n (x) \right] dx \\
 & = \int_{-1}^1 Q_n(x) \frac{d}{dx} \left[(x+1)(x^2-1) \frac{d}{dx} Q_j(x)\right] dx \\
 & = \frac{(-1)^{n+1}}{2^{n+1} (n+1)!} \int_{-1}^1 (x+1)^n (x-1)^{n+1} \frac{d^{n+2}}{dx^{n+2}} \left[(x+1)(x^2-1) \frac{d}{dx} Q_j(x)\right]  dx.
\end{align*}
A direct calculation shows 
\[
 \frac{d^{n+2}}{dx^{n+2}} \left[(x+1)(x^2-1) \frac{d}{dx} Q_j(x)\right]  = \left\{\begin{array}{ll} \frac{n(n+2)\cdot (2n+1)!}{2^{n+1} n!}, & \hbox{if} \; j=n;\\ 0, & \hbox{if}\; j<n. \end{array}\right.
\]
Thus we have $a_j = 0$ if $j < n$ and $a_n = n(n+2)$. This gives  \eqref{diff equation mod L}. Now we are ready to prove Lemma \ref{polynomial lemma even 2}. We first prove the second inequality \eqref{second inequality}. Let $P(x)$ be a polynomial of degree $\kappa$. We may write 
\[
 P(x) = \sum_{n=0}^\kappa a_n Q_n (x).
\]
We have 
\[
 \int_{-1}^{-1+\delta} (x+1)^3 |P'(x)|^2 dx  \leq \delta \int_{-1}^{-1+\delta} (x+1) (1-x^2) |P'(x)|^2 dx \leq \delta \int_{-1}^{1} (x+1) (1-x^2) |P'(x)|^2 dx.
\]
We then integrate by parts, use the differential equation and orthogonality of $\{Q_n\}$.
\begin{align*}
 \int_{-1}^{1} (x+1) (1-x^2) |P'(x)|^2 dx & = -\int_{-1}^{1} P(x) \frac{d}{dx} \left[(x+1)(1-x^2) P'(x)\right] dx\\
 & = \int_{-1}^1 \left(\sum_{n=0}^\kappa a_n Q_n(x) \right)\left(\sum_{n=0}^\kappa n(n+2) a_n (x+1) Q_n(x)\right) dx\\
 & = \sum_{n=0}^\kappa \frac{n(n+2)|a_n|^2}{2(n+1)}\\
 & \leq \kappa (\kappa+2) \int_{-1}^1 (x+1) |P(x)|^2 dx. 
\end{align*}
Combining these two inequalities, we finish the proof of \eqref{second inequality}. We then prove the first inequality \eqref{first inequality}. First of all, we have 
\begin{equation} \label{even x large} 
 \max_{x\in [0,1]} |P(x)|^2 \leq (\kappa+1)^2 \int_{0}^1 |P(x)|^2 dy \leq (\kappa+1)^2 \int_{-1}^1 (x+1) |P(x)|^2 dy.
\end{equation}
Here we apply Lemma \ref{polynomial lemma}. This deals with the case $x \in [0,1]$. Next we observe that if $x\in (-1,0)$, then we may apply a translated-version of Lemma \ref{polynomial lemma} and obtain
\[
 |P(x)|^2 \leq \max_{y \in [x,1]} |P(y)|^2 \leq \frac{(\kappa+1)^2}{1-x} \int_x^1 |P(y)|^2 dy \leq \frac{(\kappa+1)^2}{1-x^2} \int_x^1 (1+y) |P(y)|^2 dy.
\]
This immediately gives 
\[
 (1+x) |P(x)|^2 \leq \frac{(\kappa+1)^2}{1-x} \int_x^1 (1+y) |P(y)|^2 dy \leq (\kappa +1)^2 \int_{-1}^1 (1+y) |P(y)|^2 dy, \quad x \in (-1,0).
\]
Finally we combine this with the upper bound \eqref{even x large} for $x\in [0,1]$ to finish the proof of \eqref{first inequality}. 
\end{proof}

\paragraph{Decay by recursion} Finally we prove a lemma giving polynomial decay by a suitable recursion formula.   
\begin{lemma} \label{recursion lemma}
Assume that $l>1$ and $\alpha>0$ are constants. Let $S: [R,+\infty) \rightarrow [0,+\infty)$ be a function satisfying 
\begin{itemize}
 \item $S(r)\rightarrow 0$ as $r\rightarrow +\infty$;
 \item The recursion formula $S(r_2) \lesssim (r_1/r_2)^\alpha + S^l (r_1)$ holds when $r_2 \gg r_1 \gg R$.
\end{itemize}
Then given any constant $\beta \in (0, (1-1/l)\alpha)$, the decay estimate $S(r) \leq r^{-\beta}$ holds as long as $r>R_0$ is sufficiently large. 
\end{lemma}
\begin{proof}
Without loss of generality, we may assume the recursion formula 
\[
 S(r_2) \leq \frac{1}{2} (r_1/r_2)^\alpha + \frac{1}{2} S^l (r_1)
\]
holds for $r_2 \gg r_1 \gg r$. Otherwise we may slightly reduce the values of $l$ and $\alpha$. We first find a small constant $\gamma>0$ so that $S(r) \leq r^{-\gamma}$ for large $r$, then plug this estimate back in the recursion formula and slightly enlarge the value of $\gamma$, finally iterate our argument to finish the proof. We start by recalling the assumption on the limit of $S(r)$ at the infinity and choosing a large constant $M > R$ so that 
\[
 S(r) < 1/2, \qquad \forall r\in [M, M^l].
\]
This implies that we may choose a sufficiently small constant $\gamma \in (0, (1-1/l)\alpha)$ so that 
\[
 S(r) < r^{-\gamma}, \qquad \forall r \in [M, M^l].
\]
Next we prove that $S(r) \leq r^{-\gamma}$ holds for any $r\geq M$ by induction. It suffices to shows that this inequality holds for $r \in [M^{l^k}, M^{l^{k+1}}]$ if it holds for $r \in [M^{l^{k-1}}, M^{l^{k}}]$. In fact, if $r \in [M^{l^k}, M^{l^{k+1}}]$, then we have
\[
 S(r) \leq \frac{1}{2} (r^{1/l}/r)^\alpha + \frac{1}{2} S^l (r^{1/l}) \leq \frac{1}{2} r^{-(1-1/l)\alpha} + \frac{1}{2} r^{-\gamma} \leq r^{-\gamma}.
\]
Here we utilize induction hypothesis on $S(r^{1/l})$. Next we plug in $r_1 = r^{\alpha/(\alpha+\gamma l)}$ and $r_2 = r$ in the recursion formula, use the already known upper bound $S(r_1) \leq r_1^{-\gamma}$, for sufficiently large $r$, then obtain 
\[
 S(r) \leq \frac{1}{2} (r^{\alpha/(\alpha+\gamma l)}/r)^\alpha + \frac{1}{2} S^l (r^{\alpha/(\alpha+\gamma l)}) \leq r^{-\alpha \gamma l/(\alpha+\gamma l)}. 
\]
We may iterate this argument and conclude that 
\[
 S(r) \leq r^{-\gamma_k}, \qquad \forall r \geq r_k.
\]
Here $\gamma_k \in (0, (1-1/l) \alpha)$ are defined by the induction formula 
\begin{align*}
 &\gamma_0 = \gamma;& &\gamma_{k+1} = \frac{\alpha \gamma_k l}{\alpha + \gamma_k l}, \quad k \geq 0.& 
\end{align*}
In order to finish the proof, we only need to show $\gamma_k \rightarrow (1-1/l)\alpha$ as $k\rightarrow +\infty$. In fact, we may rewrite the induction formula in the form of 
\[
 (1-1/l)\alpha - \gamma_{k+1} = \frac{\alpha}{\alpha + \gamma_k l} \cdot \left[(1-1/l)\alpha -\gamma_k \right].
\]
Thus $\gamma_k \in (0, (1-1/l)\alpha)$ increases as $k \rightarrow +\infty$. This implies   
\[
 (1-1/l)\alpha - \gamma_{k+1} \leq \frac{\alpha}{\alpha + \gamma l} \cdot \left[(1-1/l)\alpha -\gamma_k \right] \qquad \Rightarrow \qquad (1-1/l)\alpha - \gamma_{k} \rightarrow 0^+. 
\]
\end{proof}

\section*{Acknowledgement}
The authors are financially supported by National Natural Science Foundation of China Project 12071339.


\begin{thebibliography}{99}
 \bibitem{newradiation} R. C\^{o}te, and C. Laurent. {``Concentration close to the cone for linear waves.''} \textit{arXiv preprint} 2109.08434. 
 \bibitem{channeleven} R. C\^{o}te, C.E. Kenig and W. Schlag. {``Energy partition for linear radial wave equation.''} \textit{Mathematische Annalen} 358, 3-4(2014): 573-607.
 \bibitem{tkm1} T. Duyckaerts, C.E. Kenig, and F. Merle. {``Universality of blow-up profile for small radial type II blow-up solutions of the energy-critical wave equation.''} \textit{The Journal of the European Mathematical Society} 13, Issue 3(2011): 533-599.
\bibitem{dkmnonradial} T. Duyckaerts, C. E. Kenig, and F. Merle. {``Universality of blow-up profile for small type II blow-up solutions of the energy-critical wave equation: the nonradial case''} \textit{The Journal of the European Mathematical Society} 14, Issue 5(2012): 1389-1454.
  \bibitem{se} T. Duyckaerts, C.E. Kenig, and F. Merle. {``Classification of radial solutions of the focusing, energy-critical wave equation.''} \textit{Cambridge Journal of Mathematics} 1(2013): 75-144.
 \bibitem{dkm3} T. Duyckaerts, C.E. Kenig, and F. Merle. {``Scattering profile for global solutions of the energy-critical wave equation.''} \textit{Journal of European Mathematical Society} 21 (2019): 2117-2162.
 \bibitem{oddtool} T. Duyckaerts, C. E. Kenig, and F. Merle. {``Decay estimates for nonradiative solutions of the energy-critical focusing wave equation.''} \textit{arXiv preprint} 1912.07655.
 \bibitem{exteriorW} T. Duyckaerts, C. E. Kenig, and F. Merle. {``Exterior energy bounds for the critical wave equation close to the ground state.''} \textit{arXiv preprint} 1912.07658.
 \bibitem{oddhigh} T. Duyckaerts, C. E. Kenig, and F. Merle. {``Soliton resolution for the critical wave equation with radial data in odd space dimensions.''}  arXiv preprint 1912.07664.
\bibitem{radiation1} F. G. Friedlander. {``On the radiation field of pulse solutions of the wave equation.''}  \textit{Proceeding of the Royal Society Series A} 269 (1962): 53-65.
\bibitem{radiation2} F. G. Friedlander. {``Radiation fields and hyperbolic scattering theory.''} \textit{Mathematical Proceedings of Cambridge Philosophical  Society} 88(1980): 483-515.
\bibitem{fourierappli} G. B. Folland. {``Fourier analysis and its applications.''} \textit{The Wadsworth and Brooks/Cole mathematics series}, 1992, Pacific Grove, California. 
\bibitem{strichartz} J. Ginibre, and G. Velo. {``Generalized Strichartz inequality for the wave equation.''} \textit{Journal of Functional Analysis} 133(1995): 50-68.
 \bibitem{channel5d} C. E. Kenig, A. Lawrie, B. Liu and W. Schlag. {``Relaxation of wave maps exterior to a ball to harmonic maps for all data''} \textit{Geometric and Functional Analysis} 24(2014): 610-647.
 \bibitem{channel} C. E. Kenig, A. Lawrie, B. Liu and W. Schlag. {``Channels of energy for the linear radial wave equation.''} \textit{Advances in Mathematics}  285(2015): 877-936.
 \bibitem{shenradiation} L. Li, R. Shen and L. Wei.  {``Explicit formula of radiation fields of free waves with applications on channel of energy''}, \textit{arXiv preprint} 2106.13396.
 \bibitem{3Dnonradiativedecay} L. Li, R. Shen and C. Wang. {``An inequality regarding non-radiative linear waves via a geometric method''}, \textit{arXiv preprint}
 \bibitem{Pecher} H. Pecher. {``Nonlinear small data scattering for the wave and Klein-Gordon equation.''} \textit{Mathematische Zeitschrift} 185(1984): 261-270.
\end{thebibliography}
\end{document}